\documentclass[a4paper,12pt]{amsart}
\usepackage{amssymb}
\usepackage{ifthen}
\usepackage[dvips]{graphicx}
\nonstopmode \numberwithin{equation}{section}
\usepackage{amssymb}
\usepackage{ifthen}
\usepackage{graphicx}
\usepackage{amsmath}
\usepackage[T1]{fontenc} 
\usepackage{comment}

\nonstopmode \numberwithin{equation}{section}
\theoremstyle{plain}
\newtheorem{thm}[equation]{Theorem}
\newtheorem{cor}[equation]{Corollary}
\newtheorem{lem}[equation]{Lemma}
\newtheorem{prop}{Proposition}

\newtheorem{conj}{Conjecture}

\newcommand{\norm}[1]{\left\lVert#1\right\rVert}

\theoremstyle{definition}
\newtheorem{defn}{Definition}[section]

\newtheorem{prob}{Problem}
\newtheorem{rem}{Remark}[section]

\newcounter{minutes}
\divide\time by 60
\newcounter{hours}
\multiply\time by 60
\addtocounter{minutes}{-\time}

\newcounter {own}
\def\theown {\thesection       .\arabic{own}}

\newenvironment{pf}[1][]{%
	\vskip 3mm
	\noindent
	\ifthenelse{\equal{#1}{}}%
	{{\slshape Proof. }}%
	{{\slshape #1.} }%
}%
{\qed\bigskip}

\newcounter{alphabet}

\newenvironment{Thm}[1][]{\refstepcounter{alphabet}%
	\bigskip%
	\noindent%
	{\bf Theorem \Alph{alphabet}}%
	\ifthenelse{\equal{#1}{}}{}{ (#1)}%
	{\bf .} \itshape}{\vskip 8pt}
\newenvironment{Lem}[1][]{\refstepcounter{alphabet}%
	\bigskip%
	\noindent%
	{\bf Lemma \Alph{alphabet}}%
	\ifthenelse{\equal{#1}{}}{}{ (#1)}%
	{\bf .} \itshape}{\vskip 8pt}

\def\be{\begin{equation}}
\def\ee{\end{equation}}

\newcommand{\bee}{\begin{enumerate}}
	\newcommand{\eee}{\end{enumerate}}

\newcommand{\blem}{\begin{lem}}
	\newcommand{\elem}{\end{lem}}
\newcommand{\bthm}{\begin{thm}}
	\newcommand{\ethm}{\end{thm}}
\newcommand{\bcor}{\begin{cor}}
	\newcommand{\ecor}{\end{cor}}
\newcommand{\beg}{\begin{examp}}
	\newcommand{\eeg}{\end{examp}}
\newcommand{\begs}{\begin{examples}}
	\newcommand{\eegs}{\end{examples}}

\newcommand{\bdefn}{\begin{defn}}
	\newcommand{\edefn}{\end{defn}}

\newcommand{\bprob}{\begin{prob}}
	\newcommand{\eprob}{\end{prob}}
\newcommand{\bei}{\begin{itemize}}
	\newcommand{\eei}{\end{itemize}}

\newcommand{\bcon}{\begin{conj}}
	\newcommand{\econ}{\end{conj}}
\newcommand{\bcons}{\begin{conjs}}
	\newcommand{\econs}{\end{conjs}}
\newcommand{\bprop}{\begin{prop}}
	\newcommand{\eprop}{\end{prop}}
\newcommand{\br}{\begin{rem}}
	\newcommand{\er}{\end{rem}}
\newcommand{\brs}{\begin{rems}}
	\newcommand{\ers}{\end{rems}}
\newcommand{\bo}{\begin{obser}}
	\newcommand{\eo}{\end{obser}}
\newcommand{\bos}{\begin{obsers}}
	\newcommand{\eos}{\end{obsers}}
\newcommand{\bpf}{\begin{pf}}
	\newcommand{\epf}{\end{pf}}
\newcommand{\ba}{\begin{array}}
	\newcommand{\ea}{\end{array}}
\newcommand{\beq}{\begin{eqnarray}}
\newcommand{\beqq}{\begin{eqnarray*}}
\newcommand{\eeq}{\end{eqnarray}}
\newcommand{\eeqq}{\end{eqnarray*}}

\begin{document}

\title{Landau-Bloch type theorem for elliptic and $K$-quasiregular harmonic mappings}

\author{Vasudevarao Allu}
\address{Vasudevarao Allu,
	School of Basic Science,
	Indian Institute of Technology Bhubaneswar,
	Bhubaneswar-752050, Odisha, India.}
\email{avrao@iitbbs.ac.in}

\author{Rohit Kumar}
\address{Rohit Kumar,
	School of Basic Science,
	Indian Institute of Technology Bhubaneswar,
	Bhubaneswar-752050, Odisha, India.}
\email{s21ma09004@iitbbs.ac.in}

\subjclass[{AMS} Subject Classification:]{Primary 30C50, 31A05; Secondary 32A18, 30C62, 30C45, 33E05.}

\def\thefootnote{}
\footnotetext{ {\tiny File:~\jobname.tex,
		printed: \number\year-\number\month-\number\day,
		\thehours.\ifnum\theminutes<10{0}\fi\theminutes }
} \makeatletter\def\thefootnote{\@arabic\c@footnote}\makeatother

\begin{abstract}
 In this paper, we establish an improved coefficient bounds for quasiregular and elliptic harmonic mappings and 
 using these bounds we establish Landau-Bloch type theorem for $(K,K')$-elliptic and K-quasiregular harmonic 
 mappings in plane. Furthermore, we prove the coefficient estimates for $K$-quasiconformal harmonic self maps defined on the unit disk $\mathbb{D}$.
\end{abstract}

\maketitle
\pagestyle{myheadings}
\markboth{Vasudevarao Allu and  Rohit Kumar}{Landau-Bloch type theorem for  elliptic and $K$-quasiregular harmonic mappings}

\section{\textbf{Introduction}}

A continuous complex-valued function $f = u + iv $ is a complex-valued harmonic function in a domain $\Omega \subset \mathbb{C}$,  if both $u$ and $v$ are real-valued harmonic functions in $\Omega$ (see \cite{Duren-2004}). The inverse function theorem and a result of Lewy \cite{Lewy-1936} shows that a harmonic function $f$ is locally univalent ({\it{i.e.}}, one-to-one) in $\Omega$ if, and only if, the Jacobian $J_f$ is non-zero in $\Omega$, where
$$
J_f(z) = |f_z(z)|^2 -|f_{\overline{z}}(z)|^2.
$$
A harmonic function $f$ is said to be sense-preserving if $J_f(z) > 0 $ for $ z $  in  $\Omega$.
Let  $\mathcal{H}$  denote  the class of  complex-valued harmonic functions $f$ in the unit disk $\mathbb{D} = \{z \in \mathbb{C}: |z| < 1\}$, normalized by $f(0) = 0 = f_z(0)-1$. Each function $f$ in $\mathcal{H}$  can be expressed as $f = h + \overline {g}$, where $h$ and $g$ are analytic functions in $\mathbb{D}$ (see \cite{Duren-2004}). Here  $h$ and $g$ are called the analytic and co-analytic parts of $f$ respectively, and  have  power series representations
$$
h(z) = z + \sum_{n = 2}^{\infty}a_n z^n \quad \mbox{and} \quad g(z) = \sum_{n = 1}^{\infty}b_n z^n.
$$
Let $\mathcal{S_H}$ be the subclass of $\mathcal{H}$ consisting of univalent, and sense-preserving harmonic mappings on $\mathbb{D}$. Let $\mathcal{H}^0 =\{ f \in \mathcal{H} : f_{\overline{z}}(0) = 0 \}$, and $\mathcal{S_{H}}^0 = \{ f \in \mathcal{S_H} : f_{\overline{z}}(0) = 0\}$. Hence for any function $f = h + \overline{g} $ in $\mathcal{H}^0$, its analytic and co-analytic parts can be represented by
$$
h(z) = z + \sum_{n = 2}^{\infty}a_n z^n \quad \mbox{and} \quad g(z) = \sum_{n = 2}^{\infty}b_n z^n
$$
respectively. It is interesting to note that $\mathcal{S_H}$ reduces to $\mathcal{S}$, the class of normalized analytic and  univalent  functions in $\mathbb{D}$, if the co-analytic part  of functions in the class  $\mathcal{S_H}$ is zero. In 1984, Clunie and Sheil-Small {\cite{Clunie-Sheil-1984}} investigated the class $\mathcal{S_H}$, together with some  geometric subclasses.  Subsequently, the class $\mathcal{S_H}$ and its subclasses have been extensively studied by several authors
(see \cite{Bshouty-Joshi-2013, Bshouty-Lyzzaik-2011, Clunie-Sheil-1984, Kalaj-Ponnusamy-Matti-2014, Wang-Liang-2001}).\\

For a continuously differentiable function $f$, denote
$$
\Lambda_f=\max _{0 \leq \theta \leq 2 \pi}\left|f_z+e^{-2 i \theta} f_{\bar{z}}\right|=\left|f_z\right|+\left|f_{\bar{z}}\right|,$$
$$
\lambda_f=\min _{0 \leq \theta \leq 2 \pi}\left|f_z+e^{-2 i \theta} f_{\bar{z}}\right|=|| f_z|-| f_{\bar{z}}|| .
$$
It is easy to see that $\left|J_f\right|=\Lambda_f \lambda_f$.\vspace{1mm}
 A harmonic mapping $f=h+\overline{g}$ defined in the unit disc $\mathbb{D}$ can be expressed as
$$
f\left(r e^{i \theta}\right)=\sum_{n=0}^{\infty} a_n r^{n} e^{i n \theta}+\sum_{n=}^{\infty} \bar{b}_n r^{n} e^{-i n \theta},\quad 0 \leq r<1,
$$
where $$h(z)=\sum_{n=0}^{\infty} {a}_{n} z^n,   \ g(z)=\sum_{n=1}^{\infty} {b}_n z^n.
$$
The classical Landau theorem (see \cite{Landau-1926}) asserts that if $f$ is a holomorphic mapping with $f(0)=0$, $f'(0)=1$ and $|f(z)|<M$ for $z \in \mathbb{D}$, then $f$ is univalent in $\mathbb{D}_{r_0}$, and $f(\mathbb{D}_{r_0})$ contains a disc $\mathbb{D}_{\sigma_0}$, where $$r_0=\frac{1}{M+\sqrt{M^2-1}}\ \ \text{and} \ \ \sigma_0= Mr_0^2.$$
  The quantities $r_0$ and $\sigma_0$ can not be improved. The extremal function is $f_0(z)= Mz \left(\frac{1-Mz}{M-z}\right)$.\vspace{2mm}

 \ \  In 2000, Chen {\it{et al.}} \cite{Chen-Gauthier-Hengertner-2000} established the following two versions of Landau-type theorems for bounded harmonic mapping on the unit disc under a suitable restriction.

\begin{Thm}\cite{Chen-Gauthier-Hengertner-2000}
Let $f$ be a harmonic mapping of the unit disc $\mathbb{D}$ such that $f(0)=0$, $f_{\bar{z}}(0)=0, f_z(0)=1$, and $|f(z)|<M$ for $z \in \mathbb{D}$. Then, $f$ is univalent on a disc $\mathbb{D}_{\rho_0}$ with
$$
\rho_0=\frac{\pi^2}{16 m M},
$$
and $f\left(\mathbb{D}_{\rho_0}\right)$ contains a schlicht disc $\mathbb{D}_{R_0}$ with
$$
R_0=\rho_0 / 2=\frac{\pi^2}{32 m M},
$$
where $m \approx 6.85$ is the minimum of the function $(3-r^2)/(r(1-r^2))$ for $0<r<1$.
\end{Thm}

\begin{Thm}\cite{Chen-Gauthier-Hengertner-2000}
 Let $f$ be a harmonic mapping of the unit disc $\mathbb{D}$ such that $f(0)=0$, $\lambda_f(0)=1$ and $\Lambda_f(z) \leq \Lambda$ for $z \in \mathbb{D}$. Then, $f$ is univalent on a disc $\mathbb{D}_{\rho_0}$ with
$$
\rho_0=\frac{\pi}{4(1+\Lambda)},
$$
and $f\left(\mathbb{D}_{\rho_0}\right)$ contains a schlicht disc $\mathbb{D}_{R_0}$ with
$$
R_0=\frac{1}{2} \rho_0=\frac{\pi}{8(1+\Lambda)} .
$$

\end{Thm}

\ \ Theorems A and B are not sharp. In 2006, better estimates for theorem A and B were given by  Grigoyan \cite{Grigoyan-2006}. In fact, Grigoyan has proved the following lemma.

\begin{Lem}\cite{Grigoyan-2006}
  Assume that $f=h+\overline{g}$ with $h$ and $g$ analytic in the unit disc $\mathbb{D}$ and let $h(z)=\sum_{n=1}^{\infty} a_n z^n$ and $g(z)=\sum_{n=1}^{\infty} b_n z^n$ for $z \in \mathbb{D}$.\vspace{2mm}
  
(a) If $|f(z)|<M$ for $z \in \mathbb{D}$, then
$$
\left|a_n\right|,\left|b_n\right| \leq M, \quad n=1,2, \ldots .
$$
\quad(b) If $\Lambda_f(z) \leq \Lambda$ for $z \in \mathbb{D}$, then
$$
\left|a_n\right|+\left|b_n\right| \leq \frac{\Lambda}{n}, \quad n=1,2, \ldots .
$$
\end{Lem}
The following result has been proved by Grigoyan using Lemma C, which improves the estimates in Theorems A and B.

\begin{Thm}\cite{Grigoyan-2006}
Let $f$ be harmonic mapping of the unit disc $\mathbb{D}$ such that $f(0)=0 , J_f(0)=1$ and $|f(z)|<M$ for $z \in \mathbb{D}$. Then, $f$ is univalent on a disc $\mathbb{D}_{\rho_1}$ with
$$
\rho_1=1-\frac{2 \sqrt{2} M}{\sqrt{\pi+8 M^2}}
$$
and $f\left(\mathbb{D}_{\rho_1}\right)$ contains a schlicht disc with
$$
R_1=\frac{\pi}{4 M}+4 M-4 M \sqrt{1+\frac{\pi}{8 M^2}}.
$$
\end{Thm}
\begin{Thm}\cite{Grigoyan-2006}
Let $f$ be harmonic mapping of the unit disc $\mathbb{D}$ such that $f(0)=0, \lambda_f(0)=1$ and $\Lambda_f(z) \leq \Lambda$ for all $z \in \mathbb{D}$. Then, $f$ is univalent on a disc $\mathbb{D}_{\rho_2}$ with
$$
\rho_2=\rho_2(\Lambda)=\frac{1}{1+\Lambda}
$$
and $f\left(\mathbb{D}_{\rho_2}\right)$ contains a schlicht disc with
$$
R_2(\Lambda)=1-\Lambda \ln \left(1+\frac{1}{\Lambda}\right).
$$
\end{Thm}	
\begin{defn}
Let $\Omega \subset \mathbb{C}$ be a doamain. A mapping $f: \Omega \rightarrow \mathbb{C}$ is said to be absolutely continuous on lines, $A C L$ in brief, in a domain $\Omega$ if for every closed rectangle $R \subset \Omega$ with sides parallel to the axes $x$ and $y, f$ is absolutely continuous on almost every horizontal line and almost every vertical line in $R$. Such a mapping $f$ has partial derivatives $f_x$ and $f_y$ {\it{a.e.}} in $\Omega$. Moreover, we say $f \in A C L^2$ if $f \in A C L$ and its partial derivatives are locally $L^2$ integrable in $\Omega$.
\end{defn}

\begin{defn}\cite{Finn-Serrin-1958}
  A sense-preserving and continuous mapping $f$ of $\mathbb{D}$ onto $\mathbb{C}$ is said to be $(K,K')$-elliptic mapping if 
  \begin{enumerate}
      \item $f$ is $ACL^2$ in $\mathbb{D}$, $J_f \neq 0$ {\it{ a.e.}} in $\mathbb{D}$,
      \item $\exists \   K\geq 1$ and $K'\geq 0$ such that $$\norm{D_f}^2\leq K J_f +K'\ {\it{a.e.}} \ \text{in}\  \mathbb{D}.$$
  \end{enumerate}
  In particular, if $K'=0$, then a $(K,K')$-elliptic mapping becomes $K$-quasiregular mapping. For more details on elliptic mappings, we refer to (\cite{Chen-Ponnusamy-Wang-2021}, \cite{Finn-Serrin-1958}, \cite{Nirenberg-1953}).
\end{defn}

It is easy to see that $\norm{D_f(z)}= \Lambda(z)$. Therefore, a sense-preserving harmonic mapping $f$ is said to be $K$-quasiregular harmonic $(K \geq 1)$ on $\mathbb{D}$ if $\Lambda_f^2 \leq K J_f$, or equivalently  $\Lambda_f \leq K \lambda_f$ holds everywhere on $\mathbb{D}$.\vspace{2mm}

If $f:G\to G'$ is $K$- quasiregular harmonic homeomorphism then $f$ is called $K$-quasiconformal harmonic mapping. \\

For $\theta \in[0,2 \pi]$ and $z \in \mathbb{D}$. Let
$$
P\left(z, e^{i \theta}\right)=\frac{1}{2 \pi} \frac{1-|z|^2}{\left|1-z e^{-i \theta}\right|^2}
$$
be the Poisson kernel. For a mapping $F \in L^1 (\partial \mathbb{D})$, the Poisson integral of $F$ is defined by
$$
f(z)=P[F](z)=\int_0^{2 \pi } P\left(z, e^{i \theta}\right) F\left(e^{i \theta}\right) d \theta .
$$

In 2007, Partyka and Sakan \cite{Partyka-Sakan-2007} proved following result.

\begin{Thm}\cite{Partyka-Sakan-2007} Given $K \geq 1$, and let $f(z)=P[F](z)$ be a harmonic $K$-quasiconformal mapping of $\mathbb{D}$ onto itself, with the boundary function $F(t)$. If $f(0)=0$, then for a.e. $z=e^{i t} \in \partial \mathbb{D}$,
$$
\frac{2^{5\left(1-K^2\right) / 2}}{\left(K^2+K-1\right)^K} \leq\left|F^{\prime}(t)\right| \leq K^{3 K} 2^{5(K-1 / K) / 2}.
$$
\end{Thm}

In 2020, Chen and Ponnusamy \cite{Chen-Ponnusamy-2020}  proved the following lemma which characterizes  elliptic mappings.

\begin{Lem}\cite{Chen-Ponnusamy-2020}
  Let $f\in \mathrm{C}^1$ be a sense-preserving mapping. Then $f$ is a  elliptic mapping if, and only if, there exist constants $k_1 \in [0,1)$ and $k_2 \in [0,\infty)$ such that
   $$\left|f_{\bar{z}}(z)\right| \leq k_1 \left|f_{z}(z)\right| +k_2  $$
    for $z \in \mathbb{D}.$ \\
    In particular, if $ f$ is $(K,K')$-elliptic mapping then 
    $$
       k_1= \frac{K-1}{K+1},\ \  k_2=\frac{\sqrt{K'}}{1+K}.
    $$
\end{Lem}
The main aim of this paper is to establish new versions of Landau-type theorems for elliptic and $K$-quasiregular harmonic mappings and some coefficient estimates for $K$-quasiconformal harmonic mappings.
\section{Main results}
We first prove coefficient estimates for elliptic harmonic mappings and $K$-quasiregular harmonic mappings.
\begin{thm} \label{Rohi-Vasu-P1-Theorem-001}
 Let $f=h+\bar{g}$ be  $(K, K')$-elliptic harmonic mapping defined on unit disc $\mathbb{D}$, where 
 \begin{equation}\label{Rohi-Vasu-P1-eq-01}
 h(z)=\sum_{n=1}^{\infty} a_n z^n \quad \text{and} \quad g(z)=\sum_{n=1}^{\infty} b_n z^n.
 \end{equation}
If $\lambda_f(z) \leq \lambda$ for $z \in \mathbb{D}$, then
\begin{equation}\label{Rohi-Vasu-P1-eq-02}
\left|a_n\right|+\left|b_n\right| \leq \frac{K \lambda+\sqrt{K'}}{n}.
\end{equation}
\end{thm} 
 We see that for the case $K=1$ and $K'=0$, we get $$\left|a_n\right|+\left|b_n\right| \leq \frac{\lambda}{n}. $$
 
Using Theorem 2.1, we obtain following corollary.
\begin{cor}\label{Rohi-Vasu-P1-Corollary-001}
 Let $f=h+\bar{g}$ be a harmonic  $K$-quasiregular mapping defined on unit disk $\mathbb{D}$, where
$$
h(z) =\sum_{n=1}^{\infty} a_n z^n \quad
\text{and}
\quad g(z) =\sum_{n=1}^{\infty} b_n z^n .
$$
If \  $\lambda_f(z) \leq \lambda$ for $z \in \mathbb{D}$, then
$$
\left|a_n\right|+\left|b_n\right| \leq \frac{K \lambda}{n}.
$$
\end{cor}
Next, we state Landau-type theorems for elliptic harmonic mappings and $K$-quasiregular mappings. 
\begin{thm}\label{Rohi-Vasu-P1-Theorem-002}
  Let $f=h+\bar{g}$ be  $(K, K')$-elliptic harmonic mapping defined on unit disc $\mathbb{D}$ such that $f(0)=0$, $\lambda_f(0)=1$ and $\lambda_f(z) \leq \lambda$ for all $z \in \mathbb{D}$. Then $f$ is univalent on a disk $\mathbb{D}_{\rho_1}$ with
$$
\rho_1=\rho_1(\lambda)=\frac{1}{1+K \lambda+\sqrt{K'}}
$$
and $f\left(\mathbb{D}_{\rho_1}\right)$ contains a schlicht disc with radius
$$
R_1=\rho_1(\lambda)+(K \lambda+\sqrt{K'})(\rho_1(\lambda)+\ln \left(\left(K \lambda+\sqrt{K'}\right) \rho_1(\lambda)\right).
$$

\end{thm}

\begin{cor}\label{Rohi-Vasu-P1-Corollary-002}
  Let $f=h+\bar{g}$ be  $K$-quasiregular  harmonic mapping defined on unit disk $\mathbb{D}$ such that $f(0)=0$, $\lambda_f(0)=1$ and $\lambda_f(z) \leq \lambda$ for all $z \in \mathbb{D}$. Then $f$ is univalent on a disc $\mathbb{D}_{\rho_1}$ with
$$
\rho_1=\rho_1(\lambda)=\frac{1}{1+K \lambda}
$$
and $f\left(\mathbb{D}_{\rho_1}\right)$ contains a schlicht disc with radius
$$
R_1=\rho_1(\lambda)+K \lambda \left(\rho_1(\lambda)+\ln \left(K \lambda 
 \ \rho_1(\lambda)\right)\right).
$$
\end{cor}
Using Corollary \ref{Rohi-Vasu-P1-Corollary-001}, we proved the following theorem.
\begin{thm}\label{Rohi-Vasu-P1-Theorem-003}
 Let $f=h+\bar{g}$ be a harmonic K-quasiregular mapping defined on unit disc $\mathbb{D}$ such that $f(0)=0, J_f(0)=1$ and $\lambda_f{(z)} \leq \lambda$ for all $z \in \mathbb{D}$. Then $f$ is univalent on a disc $\mathbb{D}_\rho$ with $$\rho(\lambda)=\frac{1}{1+\lambda K^{3 / 2}}$$
and $f\left(\mathbb{D}_\rho\right)$ contains a schlicht disc with radius $$R(\lambda)=\frac{\rho(\lambda)}{\sqrt{K}}+K \lambda\left(\rho(\lambda)+\ln (\lambda K^{3 / 2} \rho(\lambda)\right)
.$$
\end{thm}
\vspace{1mm}
Finally, using Theorem G, we prove the coefficient bounds for $K$-quasiconformal self mappings defined on the unit disc $\mathbb{D}$.

\begin{thm}\label{Rohi-Vasu-P1-Theorem-004}
 Let $f=h+\bar{g}$ be a $K$-quasiconformal harmonic self mapping on the unit disk $\mathbb{D}$, where $h$ and $g$ are analytic in $\mathbb{D}$ and given by
$$
h(z)=\sum_{n=1}^{\infty} a_n z^n \quad \text{and}
\quad g(z)=\sum_{n=1}^{\infty} b_n z^n $$ for $z \in \mathbb{D}$.
Then,
$$
\left|a_n\right|,\left|b_n\right| \leq \frac{K^{3 K} 2^{5\left(K-\frac{1}{K}\right) / 2}}{n}.
$$
\end{thm}

\section{Proof of the main results}

\begin{pf} [{\bf Proof of Theorem   \ref{Rohi-Vasu-P1-Theorem-001}}]
Let $f=h+\bar{g}$ be $(K,K')$- elliptic harmonic mapping defined on the unit disc $\mathbb{D}$, where $h$ and $g$ are given by (\ref{Rohi-Vasu-P1-eq-01}). For fixed $r \in (0,1)$ and $z=r e^{i\theta}$ for some $\theta \in [0,2\pi]$, we have 
\begin{equation}\label{Rohi-Vasu-P1-eq-03}
\frac{\partial h}{\partial r}(r e^{i \theta})=\sum_{n=1}^{\infty} n a_n r^{n-1} e^{i n \theta} \ 
\text { and } \ \frac{\partial g}{\partial r}(r e^{i \theta})=\sum_{n=1}^{\infty} n b_n r^{n-1} e^{i n \theta}.
\end{equation}
Integrating both sides of (\ref{Rohi-Vasu-P1-eq-03}) with respect to $r$ from $ 0$ to $2\pi $ after multiplying $e^{-i n \theta}$ to the left equation and multiplying $e^{i n \theta}$ to the right equation, we obtain
$$
n a_n r^{n-1}=\frac{1}{2 \pi} \int_0^{2 \pi} \frac{\partial h( re^{i \theta}) }{\partial r} e^{-i n \theta} d \theta $$
and 
$$
n b_n r^{n-1}=\frac{1}{2 \pi} \int_0^{2 \pi} \frac{\partial g( re^{i \theta}) }{\partial r} e^{i n \theta} d \theta.  $$ 
for $n=1,2...$.
Therefore, we have
$$
 n\left|a_n\right| r^{n-1} \leq \frac{1}{2 \pi} \int_0^{2 \pi} \left| {\frac{\partial}{\partial r}}h(re^{i \theta}) \right|  d \theta  $$
 and 
 $$
  n\left|b_n\right| r^{n-1} \leq \frac{1}{2 \pi} \int_0^{2 \pi} \left| {\frac{\partial }{\partial r}}g(re^{i \theta}) \right|  d \theta.$$
  Hence, we get
  $$
 n\left(\left|a_n\right|+\left|b_n\right|\right) r^{n-1} \leq \frac{1}{2 \pi} \int_0^{2 \pi}\left(\left|\frac{\partial }{\partial r}h(re^{i \theta})\right|+\left|\frac{\partial}{\partial r} g(re^{i \theta})\right|\right) d \theta 
$$
It is easy to see that
$$
\begin{aligned}
\frac{\partial }{\partial r}h( re^{i \theta}) &=\frac{\partial h}{\partial z} \frac{\partial z}{\partial r}+\frac{\partial h}{\partial \bar{z}} \frac{\partial \bar{z}}{\partial r} \\
&=h_z e^{i \theta}+0 \quad \quad \ \  \left(\because h_{\bar{z}}=0\right)\\
&=h_z e^{i \theta}
\end{aligned}
$$
and similarly, $$
\frac{\partial}{\partial r} g( re^{i \theta})=h_z e^{i \theta}.
$$
We also see that
\begin{equation}\label{Rohi-Vasu-P1-eq-04}
\left|f_z\right|=\left|h_z\right| \text { and }\left|f_{\bar{z}}\right|=\left|g_z\right|.
\end{equation}
In view of (\ref{Rohi-Vasu-P1-eq-04}), we have

\begin{equation}\label{Rohi-Vasu-P1-Eq-001}
n\left(\left|a_n\right|+\left|b_n\right|\right) r^{n-1}  \leq \frac{1}{2 \pi} \int_0^{2 \pi}\left(\left|h_z e^{i \theta}\right|+\left|g_z e^{i \theta}\right|\right) d \theta 
=\frac{1}{2 \pi} \int_0^{2 \pi}\left(\left|f_z\right|+\left|f_{\bar{z}}\right|\right) d \theta. 
\end{equation}

 Since $\lambda_f(z) \leq \lambda$ and $f$ is sense preserving, we have
$$
\left|f_z\right| \leq \lambda+\left|f_{\bar{z}}\right|.$$
Using Lemma H, we obtain
$$
\left|f_z\right| \leq \lambda+k_1\left|f_z\right|+k_2,
$$
which implies that
$$
\left|f_z\right| \leq \frac{\lambda
+k_2}{1-k_1}.
$$
Therefore, we have 
$$
\begin{aligned}
|f_z|+|f_{\bar{z}}| & \leq |f_z|+k_1|f_z|+k_2 \\
&=\left(1+k_1\right)\left|f_z\right|+k_2 \\
& \leq\left(1+k_1\right) \frac{\lambda+k_2}{1-k_1}+k_2.
\end{aligned}
$$
In view of Lemma H, we have 
\begin{equation}\label{Rohi-Vasu-P1-Eq-002}
|f_z|+|f_{\bar{z}}| \leq K \lambda+\sqrt{K^{\prime}}.
\end{equation}

By equation (\ref{Rohi-Vasu-P1-Eq-001}) and equation (\ref{Rohi-Vasu-P1-Eq-002}), we obtain
$$
\begin{aligned}
n\left(\left|a_n\right|+\left|b_n\right|\right) r^{n-1} & \leq \frac{1}{2 \pi} \int_0^{2 \pi}\left(K \lambda+\sqrt{K^{\prime}}\right) d \theta \\
&=K \lambda+\sqrt{K^{\prime}}.
\end{aligned}
$$
Now by taking $r \rightarrow 1^{-}$, we obtain
$$
\left|a_n\right|+\left|b_n\right| \leq \frac{K \lambda+\sqrt{K^{\prime}}}{n}.
$$
This completes the theorem.
\end{pf}

\begin{proof}[{\bf Proof of Corollary   \ref{Rohi-Vasu-P1-Corollary-001}}]
Let $f$ is $K$-quasiregular harmonic mappiing, then $f$ is a $(K,0)$-elliptic harmonic mapping. By substituting $K'=0$ in (\ref{Rohi-Vasu-P1-eq-02}) we obtain $$\left|a_n\right|+\left|b_n\right| \leq \frac{K \lambda}{n}.$$
\end{proof}

\begin{proof}[{\bf Proof of Theorem   \ref{Rohi-Vasu-P1-Theorem-002}}]
Let $f=h+\bar{g}$ be $(K,K')$- elliptic harmonic mapping defined on the unit disc $\mathbb{D}$, where $h$ and $g$ are given by (\ref{Rohi-Vasu-P1-eq-01}). Then
 $$
 f(z)=\sum_{n=1}^{\infty} a_n z^n+\sum_{n=1}^{\infty} \bar{b}_n \bar{z}^n \quad \text{for} \ z \in \mathbb{D} .
$$

Fix $r \in(0,1)$ and $z_1,z_2 \in \overline{\mathbb{D}}_r$ with $z_1 \neq z_2$.
It is easy to see that
$$
f\left(z_2\right)-f\left(z_1\right)=h\left(z_2\right)-h\left(z_1\right)+\overline{g\left(z_2\right)-g\left(z_1\right)}. $$
Therefore, we have
$$
|f(z_2)-f(z_1)|=\left|\int_{[z_1, z_2]} h^{\prime}(z) d z+\overline{\left(\int_{[z_1, z_2]} g^{\prime}(z) d z\right)}\right|,
$$
where $[z_1, z_2]$ is the line segment joining $z_1$ and $z_2$. It is easy to see that

\begin{align*}
\begin{split}
\left|f\left(z_2\right)-f\left(z_1\right)\right| & \geq \left| \int_{\left[z_1, z_2\right]}\left(h^{\prime}(z)-h^{\prime}(0)+h^{\prime}(0)\right) d z+\overline{\int_{\left[z_1, z_2\right]}\left(g^{\prime}(z)-g^{\prime}(0)+g^{\prime}(0)\right) d z} \right| \\
&\geq  \left| \int_{[z_1, z_2]}\left(h^{\prime}(z)-h^{\prime}(0)\right) d z+\int_{\left[z_1, z_2\right]} h^{\prime}(0) d z+ \overline{\int_{[z_1, z_2]} g^{\prime}(0) d z} \right. \\ & \ \ \quad+ \left. \overline{\int_{\left[z_1, z_2\right]}\left(g^{\prime}(z)-g^{\prime}(0)\right) d z} \right| \\ 
&\geq \left|\int_{\left[z_1, z_2\right]} h^{\prime}(0) d z+\overline{\int_{z_1, z_2]} g^{\prime}(0) d z}\right| \\ & \ 
 \ \quad-\left|\int_{\left[z_1, z_2\right]}\left(h^{\prime}(z)-h^{\prime}(0)\right) d z  +\overline{\int_{[z_1,z_2]}\left(g^{\prime}(z)-g^{\prime}(0)\right)} d z \right|
\end{split}
\end{align*}
Further,we have
\begin{align*}
\left|f\left(z_2\right)-f\left(z_1\right)\right| & \geq \lambda_f(0)\left|z_2-z_1\right|-\left|\int_{\left[z_1, z_2\right]} 
h^{\prime}\left(z\right) d z-\int_{\left[z_1, z_1\right]} h^{\prime}(0) d z\right| \\& \ \ -\left|\int_{\left[z_1, z_2\right]} g^{\prime}(z) d z-\int_{\left[z_1, z_2\right]} g^{\prime}(0) d z\right| \\
&= \lambda_f(0)\left|z_2-z_1\right|-\left|h\left(z_2\right)-h\left(z_1\right)-h^{\prime}(0)\left(z_2-z_1\right)\right| \\&
\ \ \  -\left|g\left(z_2\right)-g\left(z_1\right)-g^{\prime}(0)\left(z_2-z_1\right)\right| \\
& =|z_2-z_1| \left(1-\left|\frac{h\left(z_2\right)-h\left(z_1\right)}{z_2-z_1}-h^{\prime}(0)\right| 
-\left|\frac{g\left(z_2\right)-g\left(z_1\right)}{z_2-z_1}-g^{\prime}
(0)\right|\right).
\end{align*}
Since $z_1, z_2 \in \mathbb{\overline{D}}_r$, we obtain 
$$
\left|\frac{h\left(z_2\right)-h\left(z_1\right)}{z_2-z_1}-h^{\prime}(0)\right| \leq \sum_{n=2}^{\infty}\left|a_n\right| n r^{n-1} 
$$
and
$$
\left|\frac{g\left(z_2\right)-g\left(z_1\right)}{z_2-z_1}-g^{\prime}(0)\right| \leq \sum_{n=2}^{\infty}\left|b_n\right| n r^{n-1}.
$$
Hence 
\begin{equation}\label{Rohi-Vasu-P1-Eq-002a}
|f(z_2)-f(z_1)| \geq \left|z_2-z_1\right|\left(1-\sum_{n=2}^{\infty}\left(\left|a_n\right|+\left|b_n\right|\right) n r^{n-1}\right).
\end{equation}
Using (\ref{Rohi-Vasu-P1-Eq-002a}) and Theorem {\ref{Rohi-Vasu-P1-Theorem-001}}, we obtain
$$
\left|\frac{f(z_2)-f(z_1)}{z_2-z_1}\right| \geq 1-\left(K \lambda+\sqrt{K'}\right) \sum_{n=2}^{\infty} r^{n-1}=1-\left(K \lambda+\sqrt{K^{\prime}}\right) \frac{r}{1-r}.
$$
Since $z_2 \neq z_1$, so $f\left(z_2\right) \neq f\left(z_1\right)$ only if
$$
1-\left(K \lambda+\sqrt{K^{\prime}}\right) \frac{r}{1-r}>0 
$$
which implies that
$$
  r<\frac{1}{1+K\lambda+\sqrt{K^{\prime}}}=\rho_1.
$$
Hence, $f$ is univalent in the disc $\mathbb{D}_{\rho_1}.$
This completes the proof of the first part of the theorem.

To prove the second part of the theorem, let $|z|=\rho_1$, by Theorem {\ref{Rohi-Vasu-P1-Theorem-001}}, we obtain
\begin{align*}
|f(z)|&=\left|\sum_{n=1}^{\infty} a_n z^n+\sum_{n=1}^{\infty} \bar{b}_n \bar{z}^n\right| \\
&=\left|a_1 z+\bar{b}_1 \bar{z}^{\cdot}+\sum_{n=2}^{\infty} a_n z^n+\sum_{n=2}^{\infty} \bar{b}_n \bar{z}^n\right| \\
& \geq\left|a_1 z+\bar{b}_1 \bar{z}\right|-\left|\sum_{n=2}^{\infty} a_n z^n+\sum_{n=2}^{\infty} \bar{b}_n \bar{z}^n\right| \\
&\geq \left|\left|a_1\right|-\left|b_1\right|\right|| z|-\sum_{n=2}^{\infty}\left(\left|a_n\right|+\left|b_n\right|\right)\left|z^n\right| \\
&=\lambda_f(0)|z|-\sum_{n=2}^{\infty}\left(\left|a_n\right|+\left|b_n\right|\right)\left|z^n\right|\\
& \geq \rho_1-(K \lambda+\sqrt{K'})\sum_{n=2}^{\infty} \frac{\rho_1^n}{n}  \\
&=\rho_1+\left(K \lambda+\sqrt{K'}\right)\left(\rho_1+\ln \left(1-\rho
_1\right)\right) := R_1.
\end{align*}

Hence $f\left(\mathbb{D}_{\rho_1}\right)$ contains a schlicht disc of radius $R_1$. This completes the proof.

\end{proof}
By taking $K'=0$ in Theorem \ref{Rohi-Vasu-P1-Theorem-002}, we obtain the proof of Corollary  \ref{Rohi-Vasu-P1-Corollary-002}.
\begin{proof}[{\bf Proof of Theorem   \ref{Rohi-Vasu-P1-Theorem-003}}]
Since $J_f(0)=1$, we have 
\begin{align*}
1=|f_z(0)|^2-|f_{\bar{z}}(0)|^2 &=(|f_z(0)|+|f_{\bar{z}}(0)|)(|f_z(0)|-|f_{\bar{z}}(0)|)\\
&\leq K\left(\lambda_f(0)\right)^2
\end{align*}
which implies that
$$
 \lambda_f(0) \geq \frac{1}{\sqrt{K}}.
$$
Fix $r\in (0,1)$ and $z_1, z_2 \in \mathbb{\overline{D}}_r$ with $z_1 \neq z_2$. Then by proof of Theorem \ref{Rohi-Vasu-P1-Theorem-002} and Corollary \ref{Rohi-Vasu-P1-Corollary-001} we obtain
\begin{align*}
\left|\frac{f(z_2)-f(z_1)}{z_2-z_1}\right| &\geq \lambda_f(0)-\sum_{n=2}^{\infty}\left(|a_n|+|b_n|\right) n r^{n-1} \\
& \geq \frac{1}{\sqrt{K}}-K \lambda \sum_{n=2}^{\infty} r^{n-1}\\
& = \frac{1}{\sqrt{K}}-\frac{K \lambda r}{1-r}.
\end{align*}
Since $z_1 \neq z_2$, so $f\left(z_2\right) \neq f\left(z_2\right)$ only if
$$
\frac{1}{\sqrt{k}}-\frac{k \lambda r}{1-r} >0, $$ 
which implies that
$$ r<\frac{1}{1+\lambda k^{3 / 2}}=\rho(\lambda).
$$
Hence, $f\left(z_2\right) \neq f\left(z_1\right)$ for $r \in(0, \rho)$. This shows that $f$ is univalent in $\mathbb{D}\rho$.\\
For $|z|=\rho$, by using the proof of the Theorem \ref{Rohi-Vasu-P1-Theorem-002}, we have
\begin{align*}
|f(z)| &\geq \lambda_f(0)|z|-\sum_{n=2}^{\infty}\left(|a_n|+|b_n|\right)|z^n| \\
 & \geq \frac{\rho}{\sqrt{K}}-\sum_{n=2}^{\infty} \frac{K \lambda}{n} \rho^n \\
 & =\frac{\rho}{\sqrt{K}}+K \lambda(\rho+\ln (1-\rho))\\
 & =\frac{\rho}{\sqrt{K}}+K \lambda\left(\rho + \ln(\lambda K^{3/2} \rho)\right)=R(\lambda).
\end{align*}
Hence, $f(\mathbb{D}_\rho)$ contains a schlicht disc of radius $R(\lambda ).$ This completes the proof.
 \end{proof}

\begin{proof}[{\bf Proof of Theorem   \ref{Rohi-Vasu-P1-Theorem-004}}]
 Let $f=h+\bar{g}$ be a $K$-quasiconformal harmonic self mapping on the unit disk $\mathbb{D}$, where $h$ and $g$ are analytic in $\mathbb{D}$ and given by
$$
h(z)=\sum_{n=1}^{\infty} a_n z^n \quad \text{and}
\quad g(z)=\sum_{n=1}^{\infty} b_n z^n.$$
 Fix $r\in (0,1)$, then for $z=r e^{i\theta}$, we have 
$$
f(r e^{i \theta})=\sum_{n=1}^{\infty} a_n r^n e^{i n \theta}+\sum_{n=1}^{\infty} \bar{b}_n r^n e^{-i n \theta}.
$$
Multiplying both the sides by $e^{-i n \theta}$ and integrating with respect to $\theta$ from 0 to $2 \pi$, we get
\begin{equation}\label{Rohi-Vasu-p1-Eq-002b}
a_n r^n=\frac{1}{2 \pi} \int_0^{2 \pi} f\left(r e^{i \theta}\right) e^{-i n \theta} d \theta.
\end{equation}

Similarly,
\begin{equation}\label{Rohi-Vasu-p1-Eq-002c}
\bar{b}_n r^n=\frac{1}{2 \pi} \int_0^{2 \pi} f\left(r e^{i \theta}\right) e^{i n \theta} d \theta.
\end{equation}
It is easy to see that an integration by part of (\ref{Rohi-Vasu-p1-Eq-002b}) gives
$$
 a_n r^n=\frac{1}{2 \pi i n} \int_0^{2 \pi} f_\theta\left(r e^{i \theta}\right) e^{-i n \theta} d \theta. 
$$
which implies that
\begin{equation}\label{Rohi-Vasu-p1-Eq-003}
|a_n| r^n \leq \frac{1}{2 \pi n} \int_0^{2 \pi}\left|f_\theta(r e^{i \theta})\right| d \theta.
\end{equation}
Similarly, we see that an integration by part of (\ref{Rohi-Vasu-p1-Eq-002c}) gives
\begin{equation}\label{Rohi-Vasu-p1-Eq-004}
|b_n| r^n \leq \frac{1}{2 \pi n} \int_0^{2 \pi}\left|f_\theta(r e^{i \theta})\right| d \theta.
\end{equation}
In 2008, Kalaj \cite{Kalaj-2008} proved that the radial limits of $f_\theta$ and $f_r$ exists almost everywhere and $$\lim _{r \to 1^{-}} f_\theta\left(r e^{i \theta}\right)=f'(\theta).$$
Therefore taking $r\to 1^-$ and by Theorem F in (\ref{Rohi-Vasu-p1-Eq-003}), we get
\begin{align*}
|a_n| &\leq \frac{1}{2 \pi} \int_0^{2 \pi}\left|f'(\theta)\right| d \theta \\
& \leq \frac{K^{3 K} 2^{5(K-1 / K) / 2}}{n}.
\end{align*}
Similarly, by letting $r\to 1^-$ and by Theorem F in (\ref{Rohi-Vasu-p1-Eq-004}), we obtain
$$
|b_n| \leq \frac{K^{3 K} 2^{5(K-1 / K) / 2}}{n}.
$$
This completes the proof.
\end{proof}

\vspace{2mm}

\noindent\textbf{Acknowledgement:}  The first author thank SERB-MATRICS and the second author thank CSIR for their support.

\end{document}